\documentclass[lettersize,journal]{IEEEtran}
\usepackage{amsmath,amsfonts}
\usepackage{caption}
\usepackage{algorithmic}
\usepackage{algorithm}
\usepackage{array}
\usepackage{float}
\usepackage{enumitem}
\usepackage{multirow}
\usepackage{amssymb}
\usepackage{balance}
\usepackage{tabularx}
\usepackage[caption=false,font=normalsize,labelfont=sf,textfont=sf]{subfig}
\usepackage{tikz}
\usetikzlibrary{positioning,chains,fit,shapes,calc}
\usepackage{textcomp}

\usepackage{stfloats}
\usepackage{booktabs} 
\usepackage{url}
\usepackage{verbatim}
\usepackage{graphicx}
\usepackage{cite}
\usepackage{adjustbox} 
\usepackage{hyperref}
\usepackage{caption}

\usepackage{tikz}
\usetikzlibrary{shapes.geometric, arrows, positioning, fit, calc}
\usepackage[margin=1in]{geometry} 

\tikzstyle{largeblock} = [rectangle, draw, fill=blue!20, text centered, rounded corners, outer sep=0pt]
\tikzstyle{smallblock} = [rectangle, draw, fill=white, text centered, rounded corners, outer sep=0pt]
\tikzstyle{container} = [draw, rectangle, inner sep=0.3cm]
\tikzstyle{arrow} = [thick,->,>=stealth]
\tikzstyle{dashedcontainer} = [draw, dashed, rectangle, inner sep=0.7cm]
\tikzstyle{dottedcontainer} = [draw, rectangle, dashed, inner sep=0.7cm]
\usetikzlibrary{arrows.meta}

\usepackage{array}   

\def\BibTeX{{\rm B\kern-.05em{\sc i\kern-.025em b}\kern-.08em
    T\kern-.1667em\lower.7ex\hbox{E}\kern-.125emX}}
\begin{document}

\title{Information Geometric Framework for Comparing Point Clouds}

\author{
Amit Vishwakarma and 
KS Subrahamanian Moosath\\
Indian Institute of Space Science and Technology, Thiruvananthapuram
}

\maketitle

\begin{abstract}
In this paper, we introduce a novel method for comparing 3D point clouds, a critical task in various machine learning applications. By interpreting point clouds as 
samples from underlying probability density functions, the statistical manifold structure is given to the space of point clouds. This manifold structure will help us to
use the information geometric tools to analyze the point clouds. Our method uses the Gaussian Mixture Model (GMM) to find the probability 
density functions and the Modified Symmetric KL divergence to measure how similar the corresponding probability density functions are. This method of comparing the point clouds takes care of the geometry of the objects represented by the point clouds. To demonstrate the effectiveness of our approach, we take up five distinct case studies:(i) comparison of basic geometric shapes, (ii) comparison of 3D human body shapes within the MP FAUST dataset, (iii) comparison of animal shapes, (iv) comparison of human and animal datasets and (v) comparison of audio signals.
\end{abstract}


\section{Introduction}
\IEEEPARstart{T}{hree}-dimensional point cloud is a highly accurate digital record of an object. Point clouds provide a flexible geometric representation suitable for application in computer graphics, photogrammetry, computer vision, the construction industry, remote sensing, etc. The extraction of meaningful information from point clouds is fundamental for many applications. In particular, applications such as registration, retrieval, and autoencoding, require comparison between two or more point clouds. Comparing point clouds presents difficulties due to their inherent
 properties and the underlying surfaces they represent. Unlike structured
data on grid, point clouds do not have a common metric,
such as the Euclidean metric. Because of this the direct comparisons between point 
clouds is difficult.  When comparing the point clouds which are sampled from 3D surface, 
the primary goal is to compare the overall surface and not the specific points.
Conventional methods for point cloud comparison, including Huasdorff distance, Chamfer distance and Earth mover’s distance, have limitations. These methods directly compare the given point clouds, rely on an unstable correspondence process, and are highly sensitive to variations in sampling, which can lead to inconsistencies and inaccuracies in the comparison results.

The aim of this work is to develop a theoretical and computational framework to compare objects represented as point clouds. This method uses information geometric tools and it takes care of the shape, structure, pattern, and overall distribution.

Our method uses the deep learning models such as DGCNN (Dynamic
Graph Convolutional Neural Network) and FCNN(Fully connected Neural Network)\cite{lee2022statistical}
followed by the Gaussian Mixture Model (GMM) for giving the statistical manifold structure to the space of point clouds. Then, the similarity of the point clouds is measured using the Modified Symmetric KL divergence. In fact one can use any suitable divergence for comparing point clouds. Compared to other techniques this method provides a better understanding of point cloud similarity, considering their geometric properties. Using information geometric technique this method is effectively able to handle complex data and enhances the ability to capture intrinsic geometric properties of the point cloud data. 

To demonstrate the effectiveness of our approach, we present case studies (i) comparison of basic geometric shapes, (ii) comparison of 3D human body shapes within the MP FAUST dataset, (iii) comparison of animal shapes, (iv) comparison of human and animal datasets and (v) comparison of audio signals. The inclusion of audio signals, converted into frequency-time representations for analysis, showcases the adaptability of the information geometric method to diverse datasets.
\section{Related Works}
A significant part of working with point clouds involves comparing them to assess similarity. Commonly, the Hausdorff distance \cite{Hausdorff2008}, Chamfer distance \cite{Yang2018}, and Earth mover’s distance \cite{rubner2000earth} measures are employed. However, these measures focus primarily on how close the points are to each other, which may fail to consider the broader context and thus, the overall shape and structure of the point clouds.

The Hausdorff distance is a critical metric in the areas such as computer vision, pattern recognition, and 3D shape analysis. It is particularly useful in highlighting the dissimilarity between two point sets and is defined for two finite point sets \(P\) and \(Q\) as:

\begin{equation}
D_H(P, Q) = \max\left\{\sup_{p \in P} \inf_{q \in Q} \| p - q \|, \sup_{q \in Q} \inf_{p \in P} \| p - q \|\right\}
\end{equation}

This measure emphasizes the largest difference between two sets of 3D points, which makes it particularly sensitive to unusual points that do not fit in the pattern.

Our method leverages information geometry to mitigate the issue of unusual points by focusing on the underlying statistical properties of point clouds. By treating point clouds as samples from underlying probability distributions, we enable a comprehensive evaluation that incorporates both spatial and statistical relationships, offering a robust comparison technique.

The Chamfer distance provides an alternative that averages point-to-point distances, thereby reducing the impact of outliers. It is widely used in 3D reconstruction and mesh generation.

\begin{equation}
CD(P, Q) = \frac{1}{|P|} \sum_{p \in P} \min_{q \in Q} \| p - q \| + \frac{1}{|Q|} \sum_{q \in Q} \min_{p \in P} \| p - q \|
\end{equation}

This method calculates the average of the squared distances between each point in one set to its nearest point in the other set, offering a balance between sensitivity to large discrepancies and outliers.

By integrating geometric tools with probabilistic analysis, our method provides a nuanced assessment of point cloud similarity. Information geometric techniques enable the detection of subtle differences in point cloud shape and distribution, offering advantages in complex scenarios where context and relationships of 3D points significantly influence data structure and interpretation.

Another method for comparing the similarity and dissimilarity between point clouds is the Earth Mover’s Distance (EMD),

\begin{equation}
EMD(A, B) = \min_{\phi:A \leftrightarrow B} \sum_{a \in A} \| a - \phi(a) \|
\end{equation}

Here, \( \phi \) denotes the bijection between point clouds \(A\) and \(B\). Although EMD is insightful for measuring the dissimilarity by conceptualizing the minimum cost of transforming one point cloud into another, it might not be taking into account for the overall shape and intrinsic geometric structures of the point clouds. This limitation can lead to a less comprehensive understanding of point cloud similarities and differences, especially in contexts where the geometric and topological properties are of significant importance \cite{rubner2000earth}.

Our approach addresses this limitation by incorporating both geometric and probabilistic properties of point clouds into the comparison. By doing so, it allows for a deeper analysis of the point clouds' intrinsic structures, making our method highly effective in applications where the data's geometric and statistical characteristics are crucial. This provides a more robust comparison mechanism, especially beneficial in scenarios where understanding the complex interplay between shape, structure, and distribution is paramount.

In \cite{Yang2018},\cite{Hausdorff2008} comparing two objects reduced to the comparison between the corresponding interpoint distance matrices in the multidimensional scaling. In these methods, information about rigid similarity is only attained. But for isometric objects, allowing bend and not just rigid transformation Facundo Mémoli and Guillermo Sapiro
 developed a method in \cite{MemoliSapiro2004ComparingPC}. The underlying theory in the isometric invariant recognition is based on the Gromov-Hausdorff distance. This theory is embedded in a probabilistic setting by points sampled uniformly and using the metric as pairwise geodesic distance.
\begin{equation}
d_{GH}(X, Y) = \inf_{Z, f, g} d_{H}^{Z}(X, Y)
\end{equation}
where
\begin{equation}
d_{H}^{Z}(X, Y) = \max \left( \sup_{x \in X} d_{Z}(x, Y), \sup_{y \in Y} d_{Z}(y, X) \right),
\end{equation}
and \(X, Y\) are subsets of the manifold \(Z\) with the intrinsic geodesic distance function \(d_{Z}\). Since there is no efficient way to directly compute the Gromov-Hausdorff distance they introduced a metric
\begin{equation}
d_{J}(X,Y) = \underset{\pi\in \mathcal{P}_{n}}{\min} \, \underset{1 \leq i, j \leq n}{\max} \frac{1}{2} \left| d_{X}(x_i, x_j) - d_{Y}(y_{\pi i}, y_{\pi j}) \right|
\end{equation}

where $\mathcal{P}_{n}$ is the set of all parmutations of $\{1,2,...,n\}$, $n$ is the number of points in the point cloud.

This metric satisfies  $d_{GH}\leq d_{J}$.
The metric $d_{J}$ is computable and can be used to replace $d_{GH}$. 
The framework for comparing point clouds developed in \cite{MemoliSapiro2004ComparingPC} may face challenges with highly complex or high-dimensional geometric structures due to computational constraints. In this method, geodesic distances are approximated using techniques like the isomap algorithm, which relies on constructing a neighborhood graph from the nearest neighbors. While effective for simpler datasets, this approach can become computationally intensive for highly complex data, where the intrinsic geometric structures demand more nuanced distance calculations to accurately capture the manifold’s topology.

In the information geometric method the point clouds are represented by probability density functions and a statistical manifold structure is given to the space of point clouds. This statistical manifold structure enables a more comprehensive and robust comparison, as it considers both the spatial distribution and the statistical relationships between points.

\section{Information Geometry}
Information geometry is a new field that originated in the mid $1980$s. It explore the application of non-Euclidean geometry in probability theory in general and with a specific focus on it's application in statistical inference and estimation theory. This approach has proved to be highly effective in many applications, such as neaural network, machine learning etc.
Metrics and divergences between probability distributions is important in practical applications to look at similarity/dissimilarity among the given set of samples.

The basis of a statistical model resides in  a family of probability distributions, which are represented by a set of continuous parameters that constitute a parameter space.
 The local information contents and structures of the distributions induce certain geometrical properties on  the parameter space.   The pioneer work of Fisher in 1925, studying these geometrical properties, has received much attention.  The introduction of Riemannian metric in terms the Fisher information matrix by C. R. Rao in 1945 marked a significant milestone. This inspired many statisticians to study the geometric theory in the context of probability space.
\\
In this section, we give a brief account of the geometry of statistical manifold \cite{amari2000methods}, \cite{harsha2014f}. 

A topological space is a set equipped with a collection of subsets called open sets that satisfy certain axioms: (i) The empty set and the entire set are open (ii) The intersection of a finite number of open sets is open. The union of any collection of open sets is open.
\\
A topological space is said to be Hausdorff if for every pair of distinct points can be separated by non-intersecting open sets. This property ensures that points are individually distinguishable.
\\
 A topological space is second countable if it has a countable base for its topology, meaning that any open set can be formed as a union of sets from a countable collection.
\subsection{Smooth manifold}
An $n$-dimensional \textbf{topological manifold} $M$ is a second countable, Hausdorff topological space which is locally Euclidean. That is, for every point $p \in M$, there exists an open set $U \subset M$  and a homeomorphism $\phi:U \longrightarrow U'$, where $U'$ is an open subset of $\mathbb{R}^{n}$. $(U,\phi)$ is called a co-ordinate chart on $M$ around $p$ and $\phi$ is written as $\phi=(x^{i});i=1,...,n$. 

\begin{figure}[h]
\centering
\includegraphics[scale=.15]{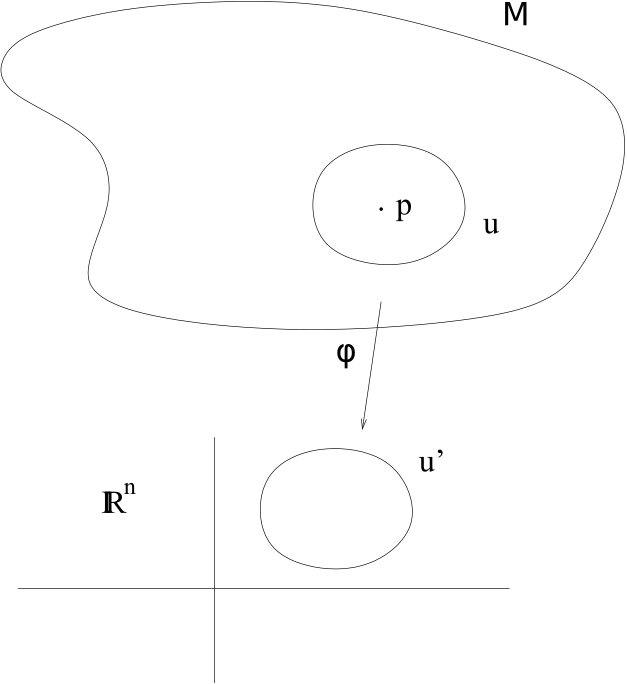}
\caption{}
\label{fig:your-label-here}
\end{figure}
If we have two charts $(U,\varphi)$ and $(V,\psi)$ on $M$ such that $U \cap V \neq \varnothing$, the composite map $\psi \circ \varphi^{-1} :\varphi(U \cap V) \longrightarrow \psi(U \cap V)$ is called the \textbf{transition map} and the charts  are said to be smoothly compatible if either $U \cap V=\varnothing$ or the transition map is a diffeomorphism. Collection of charts $\mathcal{A}$ whose domains cover $M$ is said to be an \textbf{Atlas} for M.
\begin{figure}[h]
\centering
\includegraphics[scale=.25]{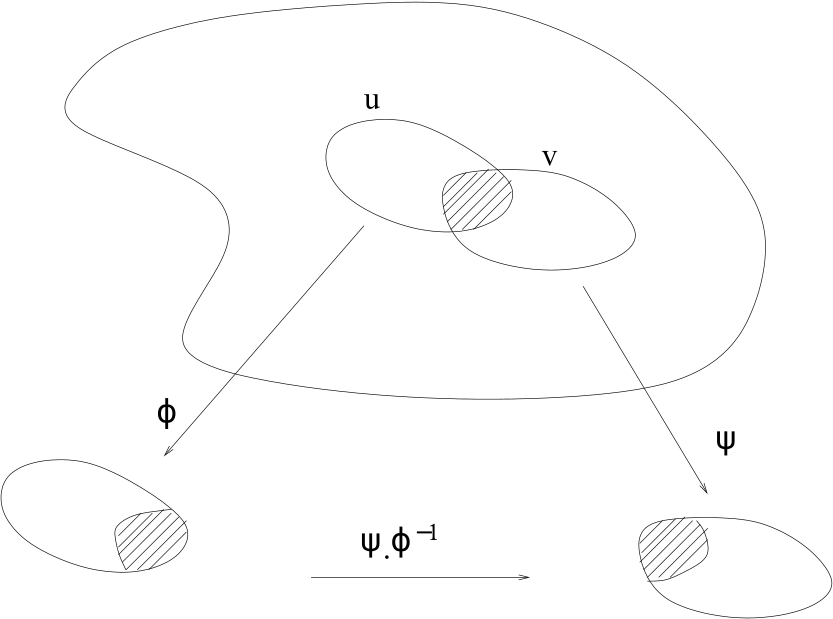}
\caption{}

\end{figure}

A smooth structure $\mathcal{A}$ is an atlas in which any two charts are smoothly compatible and $\mathcal{A}$ is maximal with respect to smooth compatibility.
  A \textbf{smooth manifold} is a pair $(M,\mathcal{A})$ where $M$ is a topological manifold and $\mathcal{A}$ is a smooth structure on $M$.

 A function $f: M \longrightarrow \mathbb{R}$  is said to be a \textbf{smooth map} if $f\circ \varphi^{-1}$ is smooth for some smooth chart $(U, \varphi)$ around each point. The set of all smooth functions from  $M$ to $\mathbb{R}$ is denoted by  $C^{\infty}(M)$ and is a real vector space.

Let $M$ be a smooth manifold and let $p \in M$. The \textbf{tangent space} $T_{p}(M)$ to $M$ at a point $p$,  is defined as the set of all derivations of $C^{\infty}(M)$ at $p$, where derivation is a linear map from $C^{\infty}(M)$ to $\mathbb{R}$ satisfying the Leibnitz rule.
 Let $(U,\phi=(x^{i}))$ be a smooth chart on $M$ around $p$. Then, $T_{p}M$ is a vector space of dimension $n$  which is spanned by $\lbrace \frac{\partial}{\partial x^{i}} \vert_{p} ; i=1,...,n\rbrace$. Each element in $T_{p}M$ is called a \textbf{tangent vector} at $p$.

A \textbf{Riemannian metric} $g=<,>$ on $M$ is a smooth symmetric 2-tensor field which is positive definite at each point that is, $g_{p}:T_{p}M \times T_{p}M \longrightarrow \mathbb{R}$ is bilinear, symmetric and positive definite. A \textbf{Riemannian manifold} is a manifold equipped with a Riemannian metric \cite{lee2013introduction}.

\subsection{Statistical manifold}
Now, we describe manifold the structure and its geometry for a statistical model.

Let $(\mathcal{X}, \Sigma, p)$ be a probability space, where $\mathcal{X} \subseteq \mathbb{R}^{n}$. Consider a family $\mathcal{S}$ of probability distributions on $\mathcal{X}$. Suppose each element of $\mathcal{S}$ can be parametrized using $n$ real-valued variables $(\theta^{1},...,\theta^{n})$ so that 
\begin{equation}
\mathcal{S}=\lbrace p_{\theta}=p(x;\theta)\;/ \;\theta=(\theta^{1},...,\theta^{n}) \in \Theta \rbrace
\end{equation}
where $\Theta$ is an open subset of $\mathbb{R}^{n}$ and the mapping $\theta \mapsto p_{\theta}$ is injective. The family $\mathcal{S}$ is called an $n$-dimensional \textbf{statistical model} or a\textbf{ parametric model}. We often write as $\mathcal{S}=\lbrace p_{\theta}\rbrace.$

For a model $\mathcal{S}=\lbrace p_{\theta} \;/ \; \theta \in \Theta \rbrace$, the mapping $\varphi:\mathcal{S}\longrightarrow \mathbb{R}^{n}$ defined by $\varphi(p_{\theta})=\theta$ allows us to consider $\varphi=(\theta^{i})$ as a coordinate system for $\mathcal{S}$. Suppose there is a $\mathit{c^{\infty}}$ diffeomorphism $\psi:\Theta\longrightarrow \psi(\Theta)$, where $\psi(\Theta)$ is an open subset of $\mathbb{R}^{n}$. Then, if we use $\rho=\psi(\theta)$ instead of $\theta$ as our parameter, we obtain $\mathcal{S}=\lbrace p_{\psi^{-1}(\rho)}\;\vert \; \rho \in \psi(\Theta) \rbrace$. This expresses the same family of probability distributions $\mathcal{S}=\lbrace p_{\theta} \rbrace$. If parametrizations which are $\mathit{c^{\infty}}$ diffeomorphic to each other is considered to be equivalent then $\mathcal{S}$ is a $\mathit{c^{\infty}}$ differentiable manifold, called the \textbf{statistical manifold}.

The tangent space $T_{\theta}(\mathcal{S})$  to the statistical manifold  $ \mathcal{S} $ at a point $p_{\theta}$ is spanned by 
\begin{equation}
\lbrace \frac{\partial}{\partial \theta^{i}} \vert_{p_{\theta}} ; i=1,...,n\rbrace.
\end{equation}
The score functions $ \lbrace \frac{\partial}{\partial \theta_{i}} \log p(x;\theta) ;i=1,...,n\rbrace$ are assumed to be linearly independent functions in $x$. Denote $\frac{\partial}{\partial \theta_{i}} \log p(x; \theta)$
 by $ \partial_{i} \ell(x, \theta) $.
The vector space $T_{\theta}^{1}(S)$ spanned by $\{\partial_{i}l(x,\theta): i=1,..n\}$ is isomorphic with $T_{\theta}(S).$ Note that $T_{\theta}(S)$ is the differential operator representation of the tangent space  and $T_{\theta}^{1}(S)$ the random variable representation of the tangent space.
\\
\noindent The expectation $E_{\theta}[\partial_{i} \ell ({x;\theta})] = \int\partial_{i} \ell({x;\theta}) p(x;\theta)dx=0$ since $p(x;\theta)$ satisfies $\int p(x;\theta)dx=1$ and so $E_{\theta}[A(x)]=0$ for  $A(x) \in T_{\theta}^{1}(\mathcal{S})$.
This expectation induces an inner product on $\mathcal{S}$ in a  natural way as $<A,B>_{\theta}\;=\;E_{\theta}[A(x)B(x)].$
Denote $\frac{\partial}{\partial \theta_{i}}$ by $\partial_{i}$
the inner product of the basis vectors $\partial_{i}$ and $\partial_{j}$ is 
\begin{eqnarray}
g_{ij}(\theta) &=& \langle \partial_{i},\partial_{j} \rangle_{\theta} \nonumber \\
&=& E_{\theta}[\partial_{i} \ell (x;\theta) \partial_{j} \ell (x;\theta)] \nonumber \\
&=& \int \partial_{i} \ell(x;\theta) \partial_{j} \ell(x;\theta) p(x;\theta) dx \label{eq:k1}.
\end{eqnarray}
 Note that the matrix $G(\theta)=(g_{ij}(\theta))$ is symmetric, for the vector $c=[c^{1},...,c^{n}]^{t}$\\
$c^{t}G(\theta)c=\int \lbrace \sum_{i=1}^{n} c^{i} \partial_{i}\ell(x;\theta)\rbrace ^2 p(x;\theta) dx \geq 0.$ 
\\
Since $\lbrace \partial_{1}\ell(x;\theta),...,\partial_{n} \ell(x;\theta) \rbrace $  are linearly independent, $G$ is positive definite. 
Hence $g=<,>$ defined above is a Riemannian metric on the statistical manifold $\mathcal{S}$,  called the \textbf{Fisher information metric}.

\subsection{Divergence Measures on Statistical Manifolds}
Divergence is a distance-like measure between two points (probability density functions) on a statistical manifold.
 The divergence $D$ on $S$ is defined as $D=D(.||.):S \times S \to \mathbb{R}$  a smooth function satisfying, for any $p,q\in S$
$$D(p||q) \geq 0\text{ and }D(p||q)=0\text{ iff }p=q.$$

The KL divergence is defined as \cite{kullback1951information},
\begin{equation}
    D_{KL}(p || q) = \int p(x; \theta_1) \log\frac{p(x; \theta_1)}{q(x; \theta_2)} dx
\end{equation}
where $p(x; \theta_1)$ and $q(x; \theta_2)$ are probability density functions. KL divergence is non-symmetric.\\
In this paper, we are using the modified version of KL divergence called
 Modified Symmetric KL divergence, denoted as \( D_{MSKL}(p \parallel q) \), is defined as

\begin{equation}
    D_{MSKL}(p \parallel q) = \frac{1}{2}\left[D_{KL}(\sqrt{p} \parallel \sqrt{q}) + D_{KL}(\sqrt{q} \parallel \sqrt{p})\right]
\end{equation}

\section{Information Geometric Framework for Point Clouds}
In this section, we establish a statistical manifold framework for point clouds. Consider the point cloud  $X$ with $n$ points in $\mathbb{R}^{m}$, i.e $X=\{x_{1},....,x_{n}:x_{i}\in \mathbb{R}^{m}\}$, assuming that all the points are distinct. Denote the set of all point clouds by $Z.$

\subsection{Covering in the Context of Point Clouds}
The concept of 'covering' plays an important role in the study of point clouds within the framework of topological spaces. A covering in the context of point clouds refers to using a group of smaller sets to represent different parts of the whole point cloud. \\
 A covering \( \mathcal{C} \) of \( X \) is defined as a finite collection of subsets \( \{U_1, U_2, \ldots, U_k\} \) of \( X \), such that $X \subseteq \bigcup_{i=1}^{k} U_i.$ Each \( U_i \) is a subset of \( X \) representing a localized region of the point cloud. The key properties to be satisfied by this covering are
\begin{enumerate}
    \item \textbf{Locality}: For each subset \( U_i \subseteq X \), there exists an open set \( N_i \) containing $U_i$ in \( X \). Mathematically, this can be represented as \( U_i \approx N_i \cap X \).

    \item \textbf{Overlap}: For any two distinct subsets \( U_i \) and \( U_j \), their intersection is not necessarily empty, that is \( U_i \cap U_j \neq \emptyset \) for some \( i \neq j \).
\end{enumerate}

\subsection{Statistical Manifold Representation:}
By choosing the finite subsets $\{U_1, U_2, \ldots, U_k\}$ of the point cloud $X$, we are simplifying the point cloud $X$ while still preserving its key features. This simplification is particularly important for finding the probability distribution governing this data using Gaussian Mixture Models (GMM). Then the statistical manifold structure is given to the space $Z$ of point clouds. By this, we obtain a manageable and accurate representation of the point clouds as a statistical manifold preserving the geometric nature of the object.

\subsection{Gaussian Mixture Model}
The fundamental idea to create the statistical manifold is to view a point cloud $X$ as a collection of $n$ samples originating from an underlying probability density function which is a point in the 
 manifold. Using the Gaussian Mixture Model \cite{jian2005robust}, a parametric probability density function $p(x; \Theta)$ is constructed for the point cloud data, where $\Theta$ is the parameter set representing the GMM.

\vspace*{5px}

 A Gaussian Mixture Model  $p(x; \Theta)$, in the context of point clouds, represents a composite distribution wherein each data point \( x \) in the point cloud \( X \) is assumed to be drawn from one of the \( K \) Gaussian components. The model is parametrized by the parameter set $\Theta = \{ (\mu_k, \Sigma_k, \pi_k): k = 1, 2, \ldots, K \}$, where $\mu_k$ is the mean vector, $\Sigma_k$ is the covariance matrix, and $\pi_k$ is the mixing coefficient.

\textit{\textbf{Mathematical Formulation:}}
\begin{enumerate}
    \item \textbf{Gaussian Distributions:} Each Gaussian component in the mixture is defined by its mean $\mu_k$ and covariance $\Sigma_k$. The probability density function of a Gaussian is given by

\begin{align}
\mathcal{N}(x|\mu_{k}, \Sigma_{k}) = &\frac{1}{\sqrt{(2\pi)^{m}|\Sigma_{k}|}} \nonumber \\
& \exp\left(-\frac{1}{2}(x - \mu_{k})^{\top} \Sigma_{k}^{-1}(x - \mu_{k})\right),
\end{align}

    where $x$ is a data point in the point cloud $X$, $m$ is the dimensionality of the data, and $|\Sigma_k|$ is the determinant of the covariance matrix.

    \item \textbf{Mixing Coefficients:} These are denoted by $\pi_k$ for each Gaussian component and they satisfy
    $
    \sum_{k=1}^{K} \pi_k = 1, \quad 0 \leq \pi_k \leq 1.
    $

    \item \textbf{Final Model:} The probability density function of the entire mixture model for a data point $x$ in the point cloud $X$ is given by
    $$
    p(x; \Theta) = \sum_{k=1}^{K} \pi_k \mathcal{N}(x|\mu_k, \Sigma_k).
    $$
\end{enumerate}
This $p(x; \Theta)$ is the statistical representation of the point cloud $X$.

\subsection{Manifold Structure for Point Clouds}

Consider the space $Z$ of point clouds, where each point cloud is represented as $X = \{x_1, \ldots, x_n: x_i \in \mathbb{R}^m\}$. Point cloud $X$ has a statistical representation $p(x,\Theta)$ where the parameter $\Theta$ varies over the parameter space $T=\{\Theta=(\mu_1,...,\mu_{k},\Sigma_{1},...,\Sigma_{K},\alpha_{1},...,\alpha_{K}):\mu_{i}\in \mathbb{R}^{m}, $$ $$ \Sigma_{i} \text{ is an } m \times m \text{ symmetric matrix, }\Sigma_{i=1}^{K}\alpha_{i}=1\}$. The set $S$ of all the probability density functions $p(x,\Theta)$ representing the point clouds in the space $Z$ is the statistical model representing $Z.$ Our aim is to give a geometric structure, called statistical manifold, to the space $Z.$

\textbf{Theorem 1:}
Let $p(x,\Theta_1) =\sum_{i=1}^{K} \alpha_{i} \mathcal{N}(x ; \mu_{i}, \sigma_i^{2})$ and $p(x,\Theta_2)=\sum_{j=1}^{L} \beta_{j} \mathcal{N}(x; \nu_{j}, \tau_j^{2})$ be two Gaussian mixture models representing the point cloud $X$ having the number of Gaussian components $K$ and $L$ respectively.  Suppose that for all $x$ in a set with no upper bound,
$$p(x,\Theta_1) = p(x,\Theta_2).$$ If the components are distinct and they are ordered lexicographically by variance and mean, then $K=L$ and for each $i$, $\alpha_{i} = \beta_{i}$, $\mu_{i} = \nu_{i}$, and $\sigma_i^{2} = \tau_i^{2}$.

\textbf{Proof:}
Arrange the components of each GMM in increasing order of their variances. For $p(x;\Theta_1)$ we have $\sigma_1^2 \leq \sigma_2^2 \leq \ldots \leq \sigma_K^2$, and for $p(x;\Theta_2)$ we have $\tau_1^2 \leq \tau_2^2 \leq \ldots \leq \tau_L^2$. If $\sigma_i^2 = \sigma_j^2$, then $\mu_i \leq \mu_j$, and similarly if $\tau_i^2 =\tau_j^2$ then $\nu_i\leq \nu_j$.

For all $x$ in a set with no upper bound,

$$\sum_{i=1}^{K} \alpha_{i} \mathcal{N}(x;\mu_{i}, \sigma_{i}^{2}) = \sum_{j=1}^{L} \beta_j \mathcal{N}(x; \nu_{j}, \tau_{j}^{2}). $$

Consider the component with the largest variance in each GMM. Without loss of generality, assume $\sigma_K^{2}$ and $\tau_L^{2}$ are the largest variances in their respective GMMs. Then,

$$ \lim_{x \to \infty} \frac{\sigma_{K}}{\mathcal{N}(x; \mu_{K}, \sigma_{K}^{2})} \sum_{i=1}^{K} \alpha_{i} \mathcal{N}(x ; \mu_{i}, \sigma_{i}^{2}) = \alpha_{K} $$
$$\lim_{x \to \infty} \frac{\tau_L}{\mathcal{N}(x | \nu_L, \tau_L^2)} \sum_{j=1}^{L} \beta_j \mathcal{N}(x | \nu_j, \tau_j^2) = \beta_L $$

Since the two GMMs are equal for all $x$, their limits must also be equal. That means $\alpha_K = \beta_L$ and $(\mu_K, \sigma_K^2) = (\nu_L, \tau_L^2)$. Then,

$$\sum_{i=1}^{K-1} \alpha_{i} \mathcal{N}(x; \mu_{i}, \sigma_{i}^{2}) = \sum_{j=1}^{L-1} \beta_{j} \mathcal{N}(x; \nu_{j}, \tau_{j}^{2})$$

By repeating the above steps we conclude that $K = L$ and for each $1 \leq i \leq K$, the parameters are equal: $\alpha_{i} = \beta_{i}$, $\mu_{i} = \nu_{i}$, and $\sigma_{i}^{2} = \tau_{i}^{2}$.
\\

\textbf{Theorem 2:}
The parameter space $$T=\{\Theta=(\mu_1,...,\mu_{k},\Sigma_{1},...,\Sigma_{K},\alpha_{1},...,\alpha_{K}):\mu_{i}\in \mathbb{R}^{m},$$ $$ \Sigma_{i} \text{ is an } m \times m \text{ symmetric matrix, }\Sigma_{i=1}^{K}\alpha_{i}=1\}$$ for GMMs with $K$ components of $m$-dimensional Gaussian is a topological manifold.\\
\textbf{Proof:} For the probability density function
$$ p(x; \Theta) = \sum_{i=1}^{K} \alpha_i \mathcal{N}(x | \mu_i, \Sigma_i),$$
each $\mu_i$ is an $m-$dimensional vector. The $m\times m$ covariance matrix $ \Sigma_i$ is symmetric, thus having $\frac{m(m+1)}{2}$ distinct elements. The total number of mixing coefficients are $K$ and since $\sum_{i=1}^{K}\alpha_{i}=1,$ only $K-1$ are independent.

The parameter space $T$ is
\begin{align*}
T = \Bigg\{ \Theta = \big( &\mu_1, \ldots, \mu_K, \Sigma_1, \ldots, \Sigma_K, \alpha_1, \ldots, \alpha_K \big) \; \Bigg| \\
&\mu_i \in \mathbb{R}^m, \Sigma_i \text{ is an } m \times m \text{ symmetric }, \nonumber \\ \text{ matrix }
&\sum_{i=1}^{K} \alpha_i = 1 \Bigg\}
\end{align*}

The total number of independent parameters in $T$ is
$ \text{dim}(T) = K \left( m + \frac{m(m+1)}{2} \right) + (K - 1).$

Note that  $T\subseteq \mathbb{R}^{dim (T)}$, so $T$ is a topological space with standard Euclidean topology. Also, around any point in $T$ there exists a neighborhood that is homeomorphic to an open subset of $\mathbb{R}^{\text{dim}(T)}$. This satisfies the local Euclidean condition for a manifold. Therefore, the parameter space $T$ of Gaussian Mixture Models with $K$ components of $m$-dimensional Gaussian is a topological manifold of dimension $K \left( m + \frac{m(m+1)}{2} \right) + (K - 1)$.
\\

\textbf{Theorem 3:}
The mapping \( h: T \rightarrow S \), defined by \( h(\Theta) = p(x; \Theta) \) for \( \Theta \in T \), is injective under the condition that the covariance matrix in the multidimensional Gaussian is diagonal.

\textbf{Proof:}
Consider parameter sets \( \Theta_1, \Theta_2 \in T \) and corresponding GMMs with their respective probability density functions 
\[ p(x; \Theta_1) = \sum_{i=1}^{K} \alpha_{i1} \mathcal{N}(x | \mu_{i1}, \Sigma_{i1}) \]
and
\[ p(x; \Theta_2) = \sum_{i=1}^{K} \alpha_{i2} \mathcal{N}(x | \mu_{i2}, \Sigma_{i2}). \]
From the Theorem 1, we know that different parameter sets lead to different GMM distributions. Extend this to the multidimensional case under the condition that the covariance matrices \( \Sigma_{i1} \) and \( \Sigma_{i2} \) are diagonal. This condition allows us to treat any linear combination of the variables in multidimensional Gaussian distributions as a univariate Gaussian distribution.

Given \( p(x; \Theta_{1}) = p(x; \Theta_{2}) \) for all \( x \), and based on the earlier Theorem 1, it implies \( \Theta_1 = \Theta_2 \). Thus, the mapping \( h: T \rightarrow S \) is injective under the condition of diagonal covariance matrix.
\\
\vspace{5px}

The parameter space $T$ is a topological manifold and the mapping $h:T\to S$ is injective hence the statistical model representing the space $Z$ of point clouds can be viewed as a statistical manifold of dimension $dim(T).$

\section{Data Overview and Preprocessing}
This section is an overview of the data sources used in the study, including basic geometrical shapes, MP FAUST dataset, and audio signals that are processed into point cloud representations. The preprocessing steps are designed to ensure that both data types are compatible with the computational framework.
\subsection{Basic Geometrical Shapes}
We consider three basic geometrical shapes in $\mathbb{R}^{3}$: unit cube $(U)$, cone $(C)$, and unit sphere $(S)$, all centered at the origin.

$U$ is defined by: $-0.5 \leq x, y, z \leq 0.5$. For the coordinates $(x_i, y_i, z_i)$ of each point $p_i$ are randomly generated using a uniform distribution over $[-0.5, 0.5]$ for each dimension.

Each coordinates $p_i$ in the cone $C$ is represented by cylindrical coordinates $(r_i,\theta_i,h_i)$, where $\theta_i$ is uniformly distributed over $[0,2\pi]$, $h$ over $[0,2]$ and the radii $r_i=1-\frac{h_i}{2}$. For $C$, each point $p_i$ is represented by its cylindrical coordinates $(r_i, \theta_i, h_i)$. The angles $\theta_i$ are uniformly distributed over $[0, 2\pi]$, the heights $h_i$ over $[0, 2]$, and the radii $r_i$ are calculated using $r_i = 1 - h_i/2$. The cylindrical coordinates are then converted into cartesian coordinates $(x_i, y_i, z_i)$ using:
$$x_i=r_{i}cos(\theta_i), y_{i}=r_{i}sin(\theta_i)\text{ and } z_i=h_i.$$

For $S$, each point $p_i$ is represented by its spherical coordinates $(\theta_i, \phi_i)$. The angles $\theta_i$ and $\phi_i$ are uniformly distributed over $[0, 2\pi]$ and $[0,\pi]$ respectively. The spherical coordinates are then converted into cartesian coordinates $(x_i, y_i, z_i)$ using:
$$x_i = sin(\phi_i)cos(\theta_i), 
y_i = sin(\phi_i)sin(\theta_i),
z_i = cos(\phi_i).$$

\subsection{Human and Animal dataset}
For the human dataset, we employ the MP FAUST dataset, a collection of high-resolution 3D scans of human figures in various poses. These point clouds capture the geometry of the human form, making them an invaluable resource for analyzing shape and structure within our geometric framework \cite{mpi_faust_dataset}.
For the animal dataset, we utilize the G-PCD: Geometry Point Cloud Dataset, which has rich local geometrical information.
\cite{epfl_geometry_point_cloud_dataset}.

\subsection{Audio Signal Processing for Point Cloud Generation}
Now present our method of generating a point cloud of an audio signal.  We start with Fourier Transform which decomposes a time-domain signal into corresponding frequencies. Then, use a Short-Time Fourier Transform (STFT) to convert audio signals into a time-frequency representation. STFT divides the audio signal into overlapping segments and then applies the Fourier Transform to each segment. This approach involves transforming the STFT data into a point cloud representation by mapping the time-frequency values onto a three-dimensional space. In the point cloud, each point corresponds to a specific pair of time-frequency with magnitude as an additional dimension. This representation enables the application of geometrical analysis techniques. Now, a brief description of the point cloud generation is given.

\subsubsection{Fourier Transform}
The Fourier Transform is a mathematical tool used to decompose a time-domain signal into its constituent frequencies. For a continuous time-domain signal \(x(t)\) and the frequency $f$ in Hz the frequency domain representation of $x(t) $ is the Fourier Transform \(X(f)\), defined as \cite{Osgood2007},
\begin{equation}
    X(f) = \int_{-\infty}^{\infty} x(t) e^{-j2\pi tf} \, dt.
\end{equation}
Note that \(X(f)\) indicating the amplitude and phase of different frequency components.

\subsubsection{Short-Time Fourier Transform (STFT)}
The STFT extends the Fourier Transform to analyze signals whose frequency content changes over time. It is defined as\cite{ADSP},
\begin{equation}
    X(n, k) = \sum_{m=-\infty}^{\infty} x(m) w(m-n) e^{-j2\pi \frac{k}{N_{\text{fft}}}m}
\end{equation}
where,
\begin{itemize}
\item \(m\) is the sample index in the original time-domain signal $x(t).$
    \item \(n\) is the time index of the windowed segment.
    \item \(k\) is the frequency bin index.
    \item \(w(m-n)\) is the window function applied to the signal, centered around time \(n\).
    \item \(N_{\text{fft}}\) is the number of Fast Fourier Transform(FFT) points, determining the resolution of the frequency bins.
\end{itemize}
The choice of window function and hop size significantly affects the analysis's resolution and accuracy. We have used  Hann window in our analysis which is defined as $\omega(n)=\frac{1}{2}\left(1-cos(\frac{2\pi n}{N_{\text{fft }}-1})\right)$\cite{leiberDifferentiableShortTimeFourier2023}

\subsubsection{Magnitude and Point Cloud Representation}
After computing the STFT, the magnitude spectrum is obtained to represent the signal's energy at various frequencies over time,
\begin{equation}
    M(n, k) = \log\left(1.0 \times 10^{-8} + |X(n, k)|\right)
\end{equation}
This magnitude is then used to form a point cloud, with each point represented as a triplet of frequency, time and magnitude.

\begin{figure}
\centering
\begin{tikzpicture}[
    node distance = 3mm and 10mm,
    start chain = going right,
    box/.style = {draw, minimum size=1.5cm, on chain},
    arrow/.style = {-latex}
]

\node[box] (a) {\includegraphics[width=.11\textwidth]{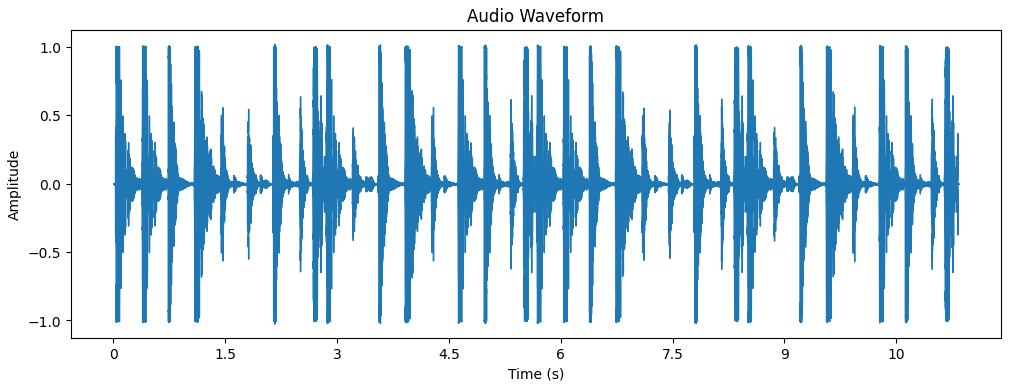}};
\node[box] (b) {\includegraphics[width=.11\textwidth]{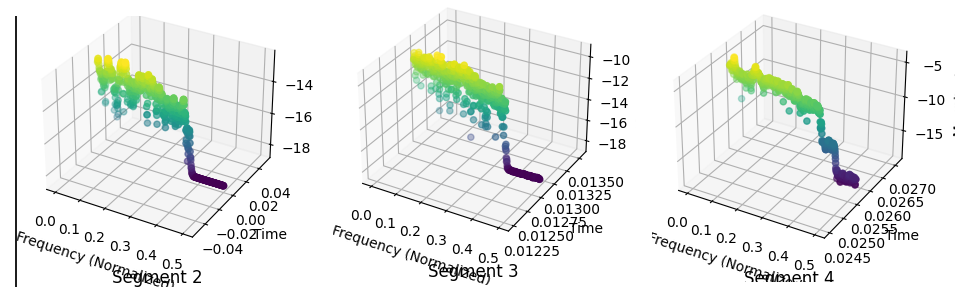}};
\node[box] (c) {\includegraphics[width=.11\textwidth]{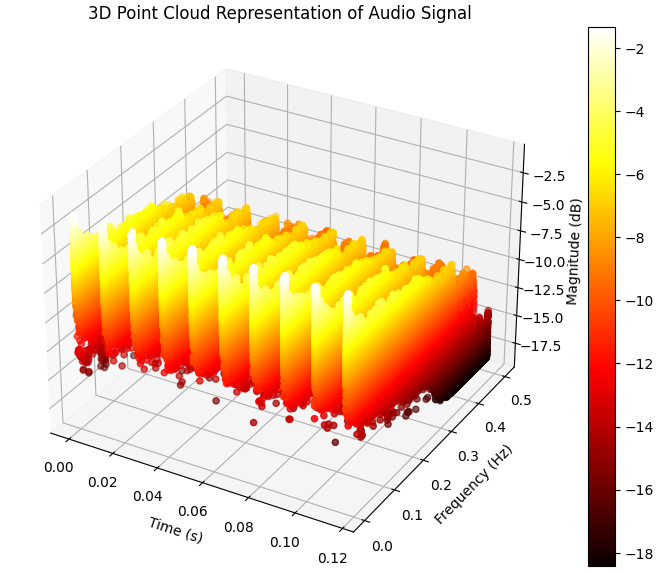}};

\draw[arrow] (a) -- (b);
\draw[arrow] (b) -- (c);

\end{tikzpicture}
\caption{Process flow from the original audio signal to the combined point cloud representation.}
\end{figure}

\subsubsection{3D Point Cloud for Each Segment}

The conversion of the magnitude spectrum into a 3D point cloud is done on a segment-by-segment basis. Each segment of the audio, traced by the STFT, is transformed into a series of 3D points where each point corresponds to a specific frequency (\(f'\)), time (\(t\)), and magnitude value (\(M(n, k)\)). This is given by,

\begin{itemize}
    \item The normalized frequency \(f' = \frac{kf_s}{N_{\text{fft}}}\), where \(k\) is the frequency bin index and \(f_s\) is the sampling rate.
    \item The time coordinate \(t = n \times \text{HopSize} / f_s\), with \(n\) being the segment index and HopSize being the step between consecutive STFT windows.
    \item The magnitude \(M(n, k)\) derived from the STFT given in equation (15).
\end{itemize}

For each STFT segment, a collection of points \((f', t, M(n, k))\) is generated, giving the spectral characteristics of that audio segment in a 3D space. This process is iterated across all segments, creating a comprehensive point cloud representation of the entire audio signal.

\subsubsection{Aggregating Point Cloud Segments into a single Representation}

After processing the audio signal through the Short-Time Fourier Transform (STFT) and converting each segment into point cloud data, next step is to aggregate these individual segments. This aggregation combines the spectral information captured in each segment to form a comprehensive point cloud representation of the entire audio signal.

(i)\hspace{0.1cm}Formation of the single Point Cloud:

The segmented point cloud data, initially organized as a tensor with dimensions \textbf{[num\_segments, num\_points\_per\_segment, 3]}, is reshaped to merge all segments. This reshaping process aligns each point from the segments into a single, cohesive data structure, effectively forming a single point cloud. Each point in this cloud is defined by three coordinates: frequency (\(f\)), time (\(t\)), and magnitude (\(M\)), corresponding to the spectral content of the audio signal at different time intervals.

(ii)\hspace{0.1cm}Saving the single Point Cloud:

For practical use, visualization, and analysis, this single point cloud is then saved in the PLY file format. The PLY format is chosen for its versatility in storing 3D data and its compatibility with a wide range of visualization and analysis tools.

\section{Computational Foundation}
In this section, the computational approach for comparing the point clouds is given. By using the deep learning techniques and EM algorithm we converted the point clouds into the probability density functions and given the statistical manifold structure to the space of point clouds. This allows to use the theory that have established in section \uppercase\expandafter{\romannumeral 4} to analyze the similarities or dissimilarities between point clouds using information geometric techniques. In the information geometric method (IGM), the similarity of the probability density function is measured using divergence, we employ the Modified symmetric KL divergence (MSKL) in this paper. In fact, one can choose a suitable divergence according to the nature of the problem.

\subsection{Sampling}

\textbf{Farthest Point Sampling (FPS) Method:}
The Farthest Point Sampling (FPS)\cite{moenning2003fast} is used to choose the samples from the point clouds. FPS is chosen because of its effectiveness in preserving the geometric and topological properties of the datasets.

Consider two distinct point clouds, $ X $ and $ Y$ of dimension $m$.
The FPS process involves,
\begin{enumerate}
    \item \textbf{Initialization}: Select an initial point $p_{1}$ randomly from the point cloud X. Set $U = \{p_{1}\}$.

    \item \textbf{Iterative Selection}: Choose the point $p_2$ that has maximum distance form $p_1$. Now choose $p_{i},$ where $i>2$, such that $p_i$ gives minimum value for $\{\max_{x\in X}|x-p_1|, \&\ldots, \ \max_{x\in X}|x-p_{i-1}|\}$.
 
Now the sample set $U$ is selected with the required number of points in this way.

    \item \textbf{Sample Extraction}: Using the above procedure choose 960 samples of 512 points each for both X and Y. Let the sample sets of X and Y be denoted by $U_{X}$ and $U_{Y}$.
    \end{enumerate}

\subsection{Labeling and Data Splitting}

In this section, the process of preparing the dataset from point clouds X and Y is given. Initially, we extracted 960 samples, each consisting of 512 points, from both the point clouds. These samples were then labeled distinctly, those from X as ``First'' and from Y as ``Second.'' This labeling was critical for identifying the source of each sample in subsequent analysis.

After the sampling combine the labeled samples into a single dataset, comprising 1920 samples in total. This dataset is then divided into training, validation, and testing sets in a 70\%, 15\%, and 15\% split, respectively. This division was carefully executed to ensure an equal distribution of samples from both point clouds in each subset, maintaining the balance and integrity of the dataset for our analysis. This setup forms the foundational basis for the computational exploration and analysis that follows.
\subsection{Model Implementation}
In this section, the implementation details of the DGCNN-FCNN (Dynamic Graph CNN followed by a Fully Convolutional Neural Network), and encoder-decoder model are discussed. This model implementation is taken from \cite{lee2022statistical}, because of its proven effectiveness in similar contexts.

This model takes input point clouds of 3$\times$512 dimension. These point clouds are passed through a series of five EdgeConv layers at the initial stage. These layers are unique in their point-wise latent space dimensions, set at $(64, 64, 128, 256)$, and a max pooling layer. We do not use the batch normalization layer. The output of this DGCNN is $1024$ dimensional feature vector. Other parameters $k=20$, leaky ReLU activation function are same as original DGCNN model.

Further, this feature vector undergoes three fully connected neural networks with dimensions $(512, 256, 2)$, with a leaky ReLU activation function and a linear output activation function. The resultant latent space is two-dimensional.

A key aspect of this approach is minimization of the reconstruction loss to ensure that the 2D latent space provides the most accurate representation of the original point cloud.

The same FCNN (fully connected neural network) is used for the decoder. The two dimensional vector from the latent space is passed through three fully connected neural networks with dimensions (256, 512, 3$\times$512), with a ReLU activation function and a linear output activation function. The final output of this process is a reconstructed three dimensional point cloud with 512 points.

Training Details:
To train our networks, we utilize the ADAM optimizer with a learning rate of 0.001. We use the Chamfer distance as the reconstruction loss. The training parameters, however, vary depending on the data type. For 3D Basic geometrical shapes, human body, and animal point cloud datasets, we conduct the training with a batch size of 16 and 400 epochs. Whereas for the audio point cloud datasets, the training is conducted over 350 epochs with a batch size of 4.

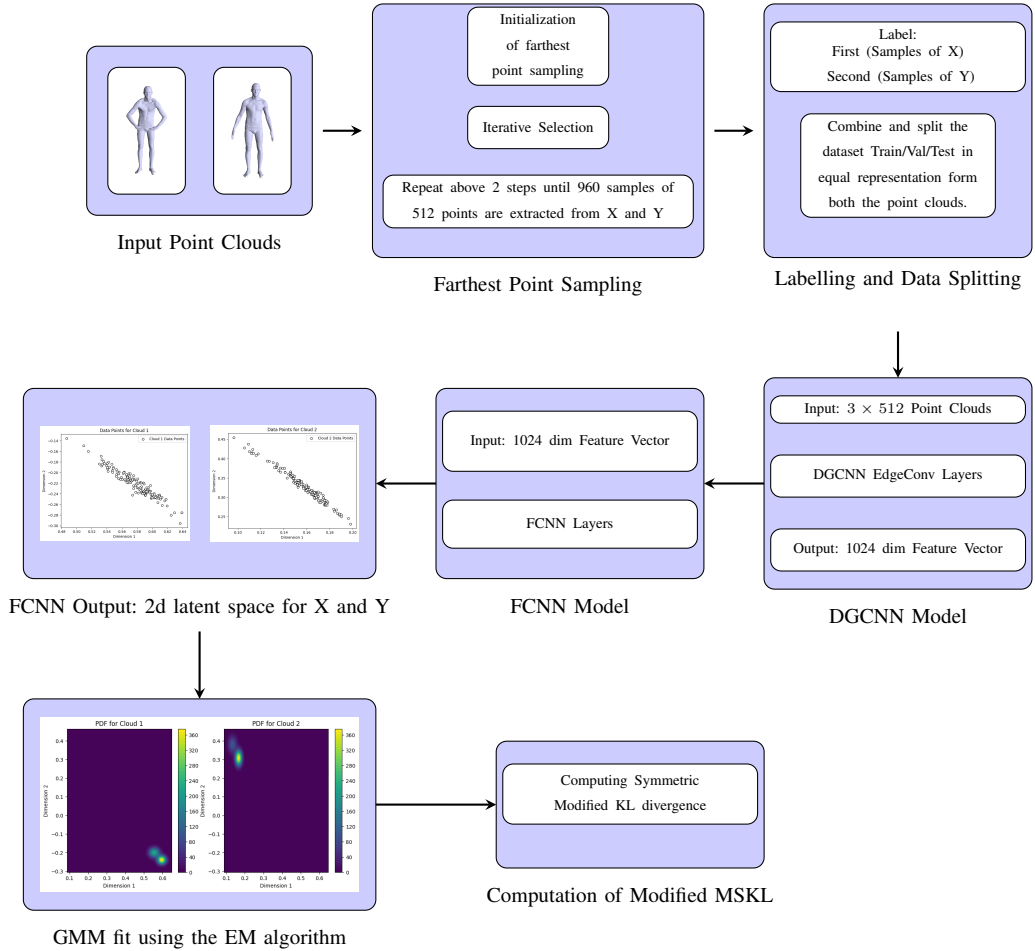
\begin{figure*}[!ht]
\centering
\hspace*{1cm}
\begin{tikzpicture}[scale=0.8, transform shape]

\node [largeblock, text width=10em, minimum height=8em] (large1) {};
\node [smallblock, text width=3em, minimum height=6em, below left=1em and 1em of large1.north east] (small1) {};
\node [smallblock, text width=3em, minimum height=6em, below right=1em and 1em of large1.north west] (small2) {};
\node at (small1) {\includegraphics[width=2em]{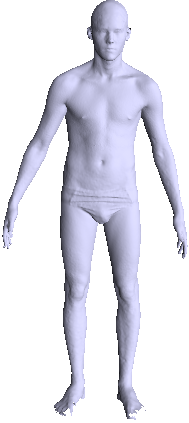}}; 
\node at (small2) {\includegraphics[width=2em]{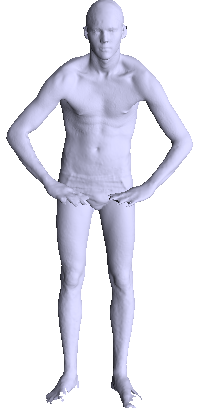}}; 

\node [largeblock, text width=15em, minimum height=12em, right=1cm of large1] (large2) {};
\node [smallblock, text width=6em, minimum height=2em] (small3) at ([yshift=-2.0em]large2.north) {\scriptsize Initialization of farthest point sampling};
\node [smallblock, text width=6em, minimum height=2em] (small4) at ([yshift=-2.0em]small3.south) {\scriptsize Iterative Selection};
\node [smallblock, text width=14em, minimum height=2em] (small5) at ([yshift=-2.5em]small4.south) {\scriptsize Repeat above 2 steps until 960 samples of 512 points are extracted from X and Y};

\draw [arrow] ([xshift=0.5em]large1.east) -- ([xshift=-0.5em]large2.west);

\node [largeblock, text width=12em, minimum height=12em, right=1cm of large2] (large3) {};
\node[smallblock, text width=4cm, minimum height=0.1cm] (small6) at([yshift=-2.5em]large3.north) {\scriptsize Label:\\ First (Samples of X)\\ \hspace{0.1cm}Second (Samples of Y)};
\node[smallblock, text width=3cm, minimum height=0.1cm] (small7) at ([yshift=-3.5em]small6.south) {\scriptsize Combine and split the dataset {Train/Val/Test} in equal representation form both the point clouds.};

\draw [arrow] ([xshift=0.5em]large2.east) -- ([xshift=-0.5em]large3.west);

\node [below=0.2cm of large1] {Input Point Clouds};
\node [below=0.2cm of large2] {Farthest Point Sampling};
\node [below=0.1cm of large3] {Labelling and Data Splitting};


\node [largeblock, text width=12em, minimum height=10em, below=2cm of large3] (large4) {};
\node[smallblock, text width =4cm, minimum height =0.2cm](small8) at ([yshift=-1.5em]large4.north){\scriptsize Input: $3\times 512$ Point Clouds};
\node[smallblock, text width =4cm, minimum height =0.7cm](small9) at ([yshift=-2.5em]small8.south){\scriptsize DGCNN EdgeConv Layers};
\node[smallblock, text width =4cm, minimum height =0.7cm](small10) at ([yshift=-2.5em]small9.south){\scriptsize Output: 1024 dim Feature Vector};

\draw [arrow] ([yshift=-3.5em]large3.south) -- (large4.north);

\node [largeblock, text width=12em, minimum height=9em, left=1cm of large4] (large5){};
\node[smallblock, text width =4cm, minimum height =1.0cm](small11) at ([yshift=-2.5em]large5.north){\scriptsize Input: 1024 dim Feature Vector};
\node[smallblock, text width =4cm, minimum height =0.8cm](small12) at ([yshift=-2.5em]small11.south){\scriptsize FCNN Layers};

\draw [arrow] (large4.west) -- (large5.east);

\node [largeblock, text width=16em, minimum height=9em, left=1cm of large5] (large6) {};

\node at ([xshift=-4em]large6.center) {\includegraphics[width=7em]{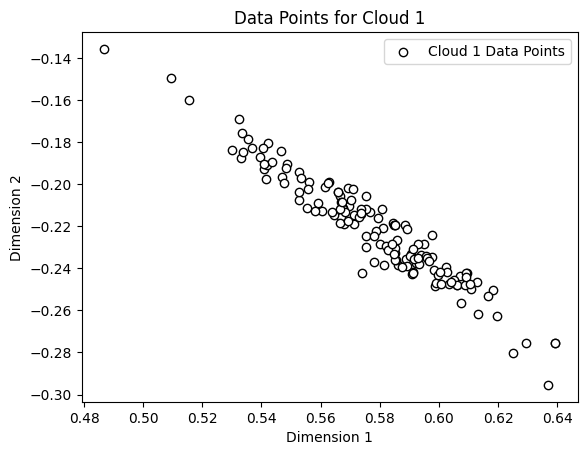}}; 

\node at ([xshift=4em]large6.center) {\includegraphics[width=7em]{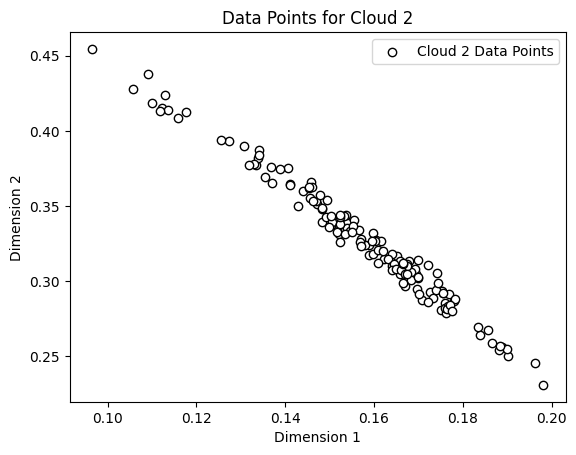}}; 

\draw [arrow] (large5.west) -- (large6.east);

\node at (large6) {};

\node [below=0.2cm of large4] {DGCNN Model};
\node [below=0.2cm of large5] {FCNN Model};
\node [below=0.2cm of large6] {FCNN Output: 2d latent space for X and Y};

\node [largeblock, text width=16em, minimum height=10em, below=2cm of large6] (large7) {};
\node at (large7.center) {\includegraphics[width=15em]{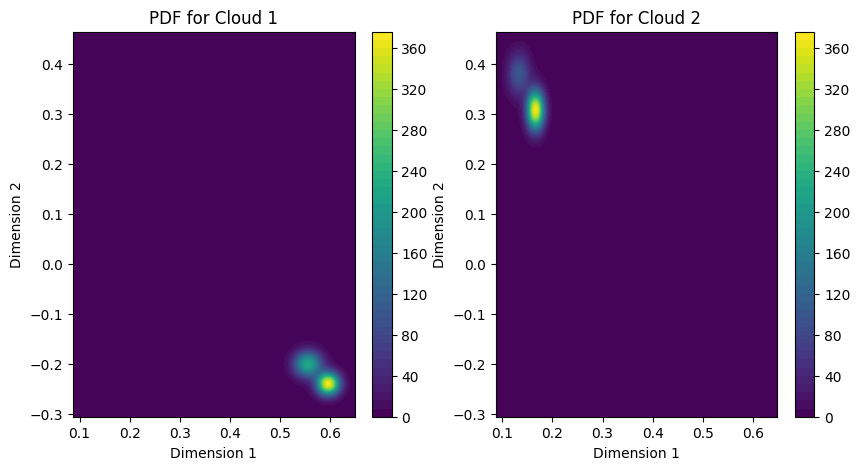}}; 
\node [below=0.2cm of large7] {GMM fit using the EM algorithm};

\draw [arrow] ([yshift=-2.5em]large6.south) -- (large7.north);

\node [largeblock, text width=12em, minimum height=6em, right=2cm of large7] (large8) {};
\node[smallblock, text width =4cm, minimum height =1.0cm](small12) at ([yshift=-2.5em]large8.north){\scriptsize Computing Symmetric Modified KL divergence};
\draw [arrow] (large7.east) -- (large8.west);
\node [below=0.2cm of large8] {Computation of Modified MSKL};

\end{tikzpicture}
\caption{Flowchart of the Computational Process.}
\label{fig:computational-process-flowchart}
\end{figure*}
\subsection{Probability Density Function Estimation Using GMM}
After the model training phase, we fit Gaussian Mixture Models (GMM) for estimating the Probability Density Functions $p(x;\Theta)$ for the point cloud data.

Post-training, test data corresponding to the label ``First'' is given as input to the trained model. This will generate a 2D latent space. Then, input the data corresponding to the label ``Second,'' to generate another 2D latent space. The Gaussian mixture model is fit to each latent space using the Expectation-maximization (EM) algorithm\cite{jannah2022parameter}.
\subsection{Space of Statistical Manifolds}

After estimating the probability density functions $p(x;\Theta)$ for the point clouds using the Gaussian mixture model, we give the statistical manifold structure to the space of point clouds using the theory that we established in Theorems 1, 2, and 3 from Section \uppercase\expandafter{\romannumeral 4}. This statistical manifold structure allows us to apply information-geometric methods to compare and analyze the point clouds in a principled and meaningful manner.

\begin{figure}[h!]
  \centering
  \includegraphics[width=0.35\textwidth]{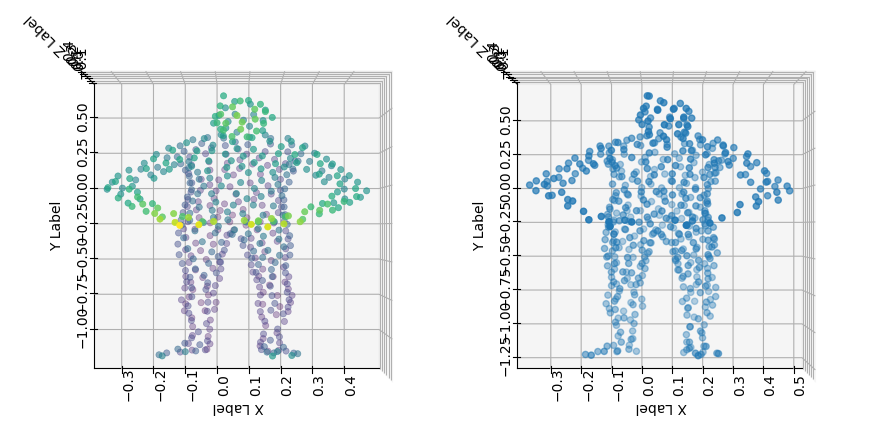}
  \caption{\small Left: The original point cloud sample extracted using FPS. Right: Reconstructed point cloud sample after passing the encoder output to the decoder model.}
  \label{fig:image1}
\end{figure}
\clearpage

\subsection{Modified Symmetric KL Divergence}
After establishing the statistical manifold structure we now use Modified Symmetric KL (MSKL) divergence for measuring similarity and dissimilarity between the point clouds. Take two sets of samples $\{x_1,.x_2,......,x_N\}$ and $\{y_1,y_2,....,y_N\}$ from two probability density functions $p(x; \Theta_1)$ and $q(y; \Theta_2)$ using a grid-based approach, instead of random sampling. This method involves creating a grid over the data range and calculating the PDF values at each grid point.

This sampling method evenly covers the entire sample area and accurately calculates the PDF values, whereas in the random sampling using Monte Carlo, it might be possible that the sampling does not cover the sample area accurately.
Then computed,
$\log\left(\frac{\sqrt{p(x_i)}}{\sqrt{q(x_i)}}\right)$ for each $x_i$ in $\{x_1, \ldots, x_N\}$ and
$\log\left(\frac{\sqrt{q(y_i)}}{\sqrt{p(y_i)}}\right)$ for each $y_i$ in $\{y_1, \ldots, y_N\}$.

  The Modified Symmetric KL divergence,
\begin{equation}
\begin{aligned}
    D_{MSKL}(p \| q) = \frac{1}{2} \Bigg[ &\sum_{i=1}^{N} \sqrt{p(x_i)}\log\left(\frac{\sqrt{p(x_i)}}{\sqrt{q(x_i)}}\right) \\
    &+ \sum_{i=1}^{N} \sqrt{q(x_i)}\log\left(\frac{\sqrt{q(y_i)}}{\sqrt{p(y_i)}}\right) \Bigg]
\end{aligned}
\end{equation}
is then computed.

\section{Case Studies}
The information geometric method is applied to five sets of datasets and analyzes the results. There are five cases each focusing on different types of datasets. Each case demonstrates the effectiveness of the information geometric method for the point cloud comparison.
\subsection{Comparision of Basic Geometrical Shapes}
We demonstrate the effectiveness of the information geometric method against the multidimensional scaling techniques (Hausdorff and Chamfer distances) for three different 3D basic shapes: sphere (S), cone (C), and cube (U). (See Table~\ref{tab:shapes}).

\begin{itemize}

\item \textbf{Accurate Identification of Identical Shapes}: The information geometric method and the multidimensional scaling techniques show the same accuracy in the case of identical shapes (e.g., sphere vs. sphere, cone vs. cone, cube vs. cube), with all values at zero.

\item \textbf{Sensitivity to Shape Variations}: The information geometric method accurately captures the shape variations. The values 6.012, 6.74 and 5.75 for sphere vs. cone, sphere vs. cube and cone vs. cube respectively shows that the point clouds are different.

The multidimensional scaling techniques show low sensitivity towards shape variation. The Hausdorff distance values for sphere vs. cone (1.41), sphere vs. cube (0.9892), and cone vs. cube (1.501) and  Chamfer distance values for the same comparisons (0.6188, 0.8675, and 1.2016) are lower than the information geometric method, suggesting a reduced sensitivity to shape variations.

\begin{table}[ht]
\centering
\caption{Comparison result of Basic Geometrical shapes}
\label{tab:shapes}
\begin{tabular}{|l|l|l|l|}
\hline
\multicolumn{1}{|c|}{\multirow{2}{*}{}} & \multicolumn{3}{c|}{\textbf{Shape}} \\ \cline{2-4} 
\multicolumn{1}{|c|}{} & \textbf{Sphere (S)} & \textbf{Cone (C)} & \textbf{Cube (U)} \\ \hline
\multirow{3}{*}{\textbf{Sphere (S)}} & Ch = 0 & & \\
 & H = 0 & & \\
 & IGM = 0 & & \\ \hline
\multirow{3}{*}{\textbf{Cone (C)}} & Ch = 0.6188 & Ch = 0 & \\
 & H = 1.41 & H = 0 & \\
 & IGM = 6.012 & IGM = 0 & \\ \hline
\multirow{3}{*}{\textbf{Cube (U)}} & Ch = 0.8675 & Ch = 1.2016 & Ch = 0 \\
 & H = 0.9892 & H = 1.501 & H = 0 \\
 & IGM = 6.74 & IGM = 5.75 & IGM = 0 \\ \hline
\end{tabular}

\end{table}
\end{itemize}
\subsection{Audio Data Analysis}
In the audio analysis, we transferred the audio signals into point clouds using the Short-Time Fourier Transform (STFT).

Below we discuss the merits of the information geometric method compared to the multidimensional scaling techniques.(See Table~\ref{tab:audio}).
\begin{itemize}
\item \textbf{Accurate Identification of Identical Audio Tracks}: The information geometric method and the multidimensional scaling techniques show the same accuracy in the case of identical audio signals (e.g., Audio1 vs. Audio1), with all values at zero.

\item \textbf{Balanced Sensitivity to Audio Variations}:  Information geometric method shows balanced sensitivity to audio variations. For different audio signals (Audio1 vs. Audio2), the value is 0.5999, indicating two audios are different. Also for the same sentence by different people, the information geometric method gives a value of 0.2239. It shows the method captures the geometry of the signal and the small variation effectively.

On the other hand, the multidimensional scaling techniques, facing difficulty in capturing audio variations. The Hausdorff distance value of 0.4698 for Audio1 vs. Audio2 is lower than the IGM value, while the Chamfer distance value of 0.0542 for the same comparison is very small shows two audios are similar which is not the case.

 \item \textbf{Differentiation of Pitch and Tempo Variations}: Information geometric method effectively captures differences in pitch and tempo. The information geometric method shows values 0.10251 (similar pitch) and 0.2667 (different pitch) when the same person speaks different sentences in different pitch.

In comparison, the multidimensional scaling techniques under-represent the impact of pitch and tempo variations. The Hausdorff distance shows values 0.2173 (similar pitch) and  0.2891 (different pitch), while the Chamfer distance shows 0.00742 (similar pitch) and 0.0096 (different pitch), showing the superiority of the information geometric method.

\end{itemize}

\begin{table}[H]
  \centering
  \caption{Comparison result of Audio signals}
  \label{tab:audio}
  \begin{tabular}{|p{0.3\linewidth}|p{0.15\linewidth}|p{0.15\linewidth}|p{0.15\linewidth}|}
    \hline
    & IGM & Hausdorff & Chamfer \\
    \hline
    Audio1 - Audio1 & 0 & 0 & 0 \\
    \hline
    Audio2 - Audio2 & 0 & 0 & 0 \\
    \hline
    Audio1 - Audio 2 & 0.5999 & 0.4698 & 0.0542 \\
    \hline
    Same sentence but two different people & 0.2239 & 0.7361 & 0.01072 \\
    \hline
    Same person two different sentence (similar pitch) & 0.10251 & 0.2891 & 0.0096 \\
    \hline
    Same person two different sentence (Different Pitch) & 0.2667 & 0.2173 & 0.00742 \\
    \hline
  \end{tabular}
\end{table}

\subsection{Point cloud comparison of the Human body:}

The analysis of point clouds of human body dataset \cite{mpi_faust_dataset} shows that the information geometric method effectively handles complex point cloud data compared to the multidimensional scaling technique. (See Table~\ref{tab:3d Body}).

\begin{itemize}

\item \textbf{Accurate Identification of Identical Shapes}: The information geometric method and the multidimensional scaling technique show equal capability when applied to the identical point clouds (e.g., Man1 vs. Man1 and Man2 vs. Man2).

\item \textbf{Balanced Sensitivity to Shape Deformations}: When comparing the same person in different postures (e.g., Man1 vs. Man2), the information geometric method shows a balanced output of 0.1412, suggesting slight yet small distinctions due to variations in body posture. In contrast, the multidimensional scaling techniques exhibit limitations in capturing shape deformations. The Hausdorff distance (0.2661) somewhat overstates the differences and the Chamfer distance (0.0888) underestimates them.

\item \textbf{Differentiation of Topologically Similar Shapes}: The information geometric method shows it's capability in differentiating topologically similar shapes (e.g., Man vs. Woman).  It accurately captures the subtle differences, with 0.2645 (Man2 vs. Woman1) which is somewhat higher than the value between Man1 vs. Man2 as expected, demonstrating its ability to capture geometric changes. On the other hand, the multidimensional scaling techniques, Hausdorff (0.4423) and Chamfer (0.2875), slightly overestimate the differences.

\end{itemize}
\subsection{Point cloud comparison of Animals:}

The analysis of the animal dataset \cite{epfl_geometry_point_cloud_dataset}  highlights the effectiveness of our information geometric method in capturing the differences between different animal point clouds. (See Table~\ref{tab:3d Body})

\begin{itemize}

\item \textbf{Accurate Identification of Identical Shapes}: The information geometric method and multidimensional scaling technique show equal capability when applied to the identical point clouds (e.g., Rabbit vs. Rabbit and Dragon vs. Dragon).

\item \textbf{Differentiation of Topologically Similar Shapes}: The information geometric method effectively distinguishes between different animal shapes (e.g., Rabbit vs. Dragon), with value of 0.4512. The value shows the significant differences between the two animal shapes while maintaining a balanced approach. In comparison, the multidimensional scaling techniques, the Chamfer distance (0.2779) was somewhat lower and the Hausdorff distance (0.4598) was almost the same as our output.

\end{itemize}
\subsection{Point cloud comparison of the Human body and Animals:}
Below we discuss the comparison of two completely different datasets: the human body dataset and the animal body dataset. (See Table~\ref{tab:3d Body}).

\begin{itemize}
\item \textbf{Differentiating Between Human and Rabbit Point Clouds}: The information geometric method and the multidimensional scaling technique effectively differentiate between human and rabbit shapes. The inforamation geometric method, Chamfer distance and Hausdorff distances are consistently high which shows the human and rabbit point clouds are different.
\item \textbf{Differentiating Between Human and Dragon Point Clouds}: The information geometric method and the multidimensional scaling techniques exhibit their effectiveness in distinguishing between human and dragon shapes. In all the cases the values are high.

\end{itemize}

\begin{table*}[tp]
\small
\begin{center}
\caption{Comparison result of complex shapes}
\label{tab:3d Body}
\begin{adjustbox}{width=\textwidth}
\begin{tabular}{|c|c|c|c|c|c|c|c|c|}
\hline
 & \textbf{Man1} & \textbf{Man2} & \textbf{Woman1} & \textbf{Rabbit} & \textbf{Dragon} \\ \hline
 & \includegraphics[width=0.8cm]{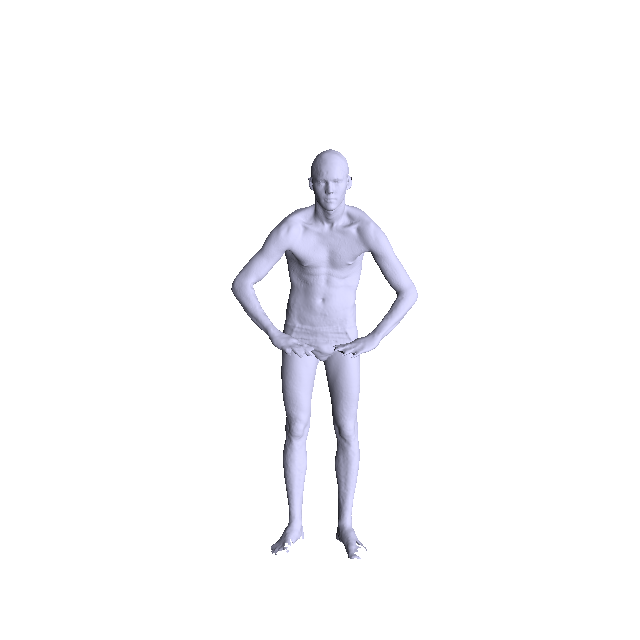} & \includegraphics[width=0.8cm]{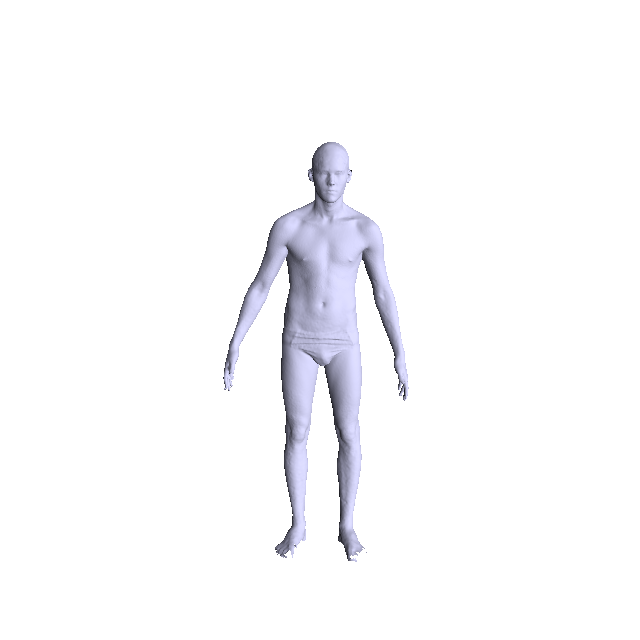} & \includegraphics[width=0.8cm]{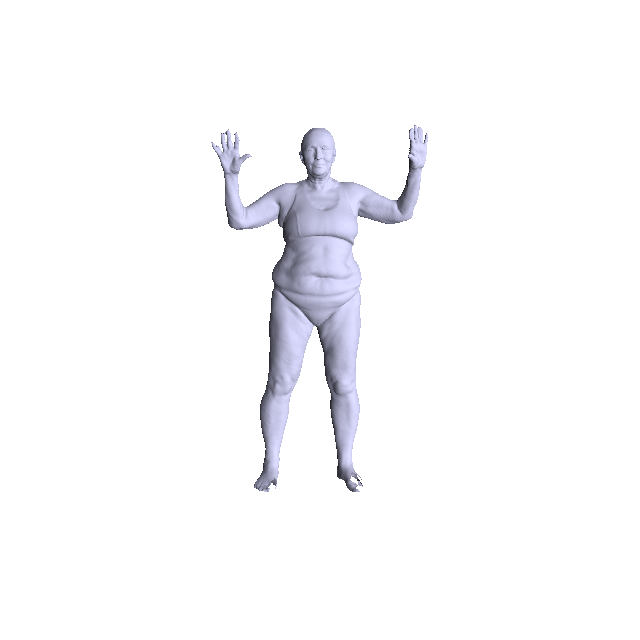} & \includegraphics[width=0.8cm]{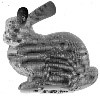} & \includegraphics[width=0.8cm]{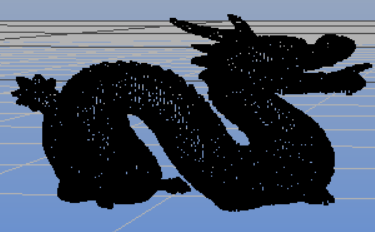} \\ \hline

\includegraphics[width=0.8cm]{tr_scan_001.png} & 
\begin{tabular}{@{}c@{}}
   \scriptsize Ch= 0 \\
\scriptsize H= 0 \\
\scriptsize IGM= 0
\end{tabular} & & & &
\\ \hline

\includegraphics[width=0.8cm]{tr_scan_000.png} & 
\begin{tabular}{@{}c@{}}
    \scriptsize Ch=0.0888 \\
   \scriptsize H=0.2661 \\

    \scriptsize IGM=0.1412
\end{tabular} & 
\begin{tabular}{@{}c@{}}
    \scriptsize Ch= 0 \\
    \scriptsize H= 0 \\

    \scriptsize IGM= 0
\end{tabular} & & & 
\\ \hline

\includegraphics[width=0.8cm]{test_scan_162.png} & 
\begin{tabular}{@{}c@{}}
   \scriptsize Ch=0.2615 \\
   \scriptsize H=0.4563 \\

   \scriptsize IGM=0.2359
\end{tabular} & 
\begin{tabular}{@{}c@{}}
   \scriptsize Ch=0.2875 \\
   \scriptsize H=0.4423 \\

   \scriptsize IGM=0.2645
\end{tabular}
&\begin{tabular}{@{}c@{}}
   \scriptsize Ch=0 \\
   \scriptsize H=0 \\

   \scriptsize IGM =0
\end{tabular} & &
\\ \hline
\includegraphics[width=0.8cm]{raaaa.jpg} & 
\begin{tabular}{@{}c@{}}
  \scriptsize Ch= 0.7879 \\
   \scriptsize H= 1.2488 \\

 \scriptsize IGM= 0.5871
\end{tabular} & 
\begin{tabular}{@{}c@{}}
  \scriptsize Ch= 0.8094 \\
   \scriptsize H= 1.2461 \\

 \scriptsize IGM = 0.5478
\end{tabular} &

\begin{tabular}{@{}c@{}}
  \scriptsize Ch= 0.9030 \\
   \scriptsize H= 1.1394 \\

 \scriptsize IGM = 0.3612
\end{tabular} &
\begin{tabular}{@{}c@{}}
  \scriptsize Ch= 0 \\
   \scriptsize H= 0 \\
 
 \scriptsize IGM= 0
\end{tabular} &

\\ \hline
\includegraphics[width=0.8cm]{dragon.png} & 
\begin{tabular}{@{}c@{}}
   \scriptsize Ch=0.5291 \\
   \scriptsize H= 1.2073 \\

   \scriptsize IGM =0.5466
\end{tabular} & 
\begin{tabular}{@{}c@{}}
   \scriptsize Ch=0.5626 \\
   \scriptsize H=01.2015 \\
 
   \scriptsize IGM =0.5789
\end{tabular} &
\begin{tabular}{@{}c@{}}
   \scriptsize Ch=0.5994 \\
   \scriptsize H=1.0161 \\
  
   \scriptsize IGM=0.6124
\end{tabular} &
\begin{tabular}{@{}c@{}}
   \scriptsize Ch=0.2779 \\
   \scriptsize H=0.4598 \\

   \scriptsize IGM =0.4512
\end{tabular} &
\begin{tabular}{@{}c@{}}
   \scriptsize Ch=0 \\
   \scriptsize H=0 \\

   \scriptsize IGM =0
\end{tabular} \\
\hline
\end{tabular}
\end{adjustbox}
\end{center}

\label{tab: Table 1}
\end{table*}

The results highlight the overall effectiveness of our information geometric method, in accurately capturing shape deformations, differentiating topologically similar shapes, and distinguishing between completely different point clouds giving a well-balanced output in 3D basic geometrical shapes, human body and animal case. In audio data analysis, our method effectively captures audio variations, pitch, and tempo differences. However, the multidimensional scaling techniques have certain limitations in analyzing both shape and audio point cloud data. In certain instances, it has a tendency to either overstate or underestimate the disparities, hence lacking the ability to accurately capture the geometry of the data.

\section{Conclusion}
Information geometric method discussed in this paper gives a more accurate and balanced comparison of the point clouds and captures the geometry effectively. By giving the statistical manifold structure to the space of point clouds the geometric tools can be applied to datasets for comparing them efficiently. This method is adaptable to various kinds of data and able to give a geometric framework for comparing the datasets. The advantages of the information geometric method over other techniques are, (i) it offers a mathematically rigorous framework based on information geometry, (ii) it is a comprehensive method that can be used for different types of datasets and (iii) the method's balanced sensitivity and robustness to variations in geometry make it suitable for many practical applications.

The information geometric method consistently exhibits a balanced sensitivity to shape deformations and topological variations in all the different types of datasets including basic geometrical shapes, audio signals, human body scans, and animal point cloud. Moreover, the information geometric method shows robustness in handling complex datasets, such as audio signals with variations in pitch and tempo. It effectively distinguishes between audio of different speakers or sentences, showing its ability to capture complex datasets' geometry. In the context of 3D shape analysis, information geometric method is effective in differentiating topologically similar shapes, such as human body scans in different postures.

Thus, the information geometric approach presented in this paper provides a new and effective way to analyze and compare point clouds. The method effectively handles different types of datasets which could be very useful in areas like audio processing, computer vision, 3D modeling, etc. The method also suggests that one can classify objects by fixing a threshold. In particular, in the medical field, diseases like tumors, etc can be meticulously detected. \\

In future work, our plan is to take up some interesting research problems in image processing, remote sensing, etc, using the information geometric method that we have establish in this paper. For the computational cost, here we reduced the dimension of the 3D dataset to 2D latent space, in future we will try to work with the three dimensional latent space. We will also explore the concept of connection which helps in understanding how the geometry changes when we move along the manifold. So along with the divergence measure, we also like to incorporate the curvature measure so as to capture the geometry of the data more effectively.

\section{Acknowledgment}
Amit Vishwakarma is thankful to the Indian Institute of Space Science and Technology, Department of Space, Govt. of India for the award of the doctoral research fellowship.

\bibliographystyle{plain} 
\bibliography{mybib} 

\end{document}